\documentclass[aos,MSNbibl,seceqn,citesort,dvips]{arximspdf}
\usepackage{dcolumn}
\usepackage{graphicx}

\doi{10.1214/11-AOS948}
\volume{40}
\issue{1}
\pubyear{2012}
\firstpage{104}
\lastpage{130}

\makeatletter

\newcolumntype{d}[1]{D{.}{.}{#1}}

\newtheorem{Theorem}{Theorem}[section]
\newtheorem{Proposition}[Theorem]{Proposition}
\newtheorem{Corollary}[Theorem]{Corollary}

\newproclaim{Example}{Example}[section]

\newcommand{\bfY}{\mathbf{Y}}
\newcommand{\ind}{\mathbb{I}}

\newcommand{\M}{\mathcal{M}}
\newcommand{\te}{\theta}
\newcommand{\bfs}{\mathbf{s}}

\newcommand{\cS}{\mathcal{S}}
\newcommand{\bfone}{\mathbf{1}}
\newcommand{\bfZ}{\mathbf{Z}}
\newcommand{\bfV}{\mathbf{V}}
\newcommand{\bfx}{\mathbf{x}}

\newcommand{\bfU}{\mathbf{U}}
\newcommand{\bff}{\mathbf{f}}
\newcommand{\bfW}{\mathbf{W}}
\newcommand{\mbalpha}{\bolds\alpha}

\newcommand{\bfi}{\mathbf{i}}
\newcommand{\bfa}{\mathbf{a}}
\newcommand{\D}{\mathcal{D}}
\newcommand{\raw}{\rightarrow}
\newcommand{\nti}{n \to\infty}
\newcommand{\hT}{\hat{T}}
\newcommand{\hG}{\hat{G}}

\newcommand{\hW}{\hat{W}}
\newcommand{\bfy}{\mathbf{y}}
\newcommand{\hU}{\hat{U}}
\newcommand{\hth}{\hat{\theta}}
\newcommand{\bft}{\mathbf{t}}
\newcommand{\C}{\mathcal{C}}
\newcommand{\bfzero}{\mathbf{0}}
\newcommand{\N}{\mathcal{N}}
\newcommand{\rl}{\mathbb{R}}
\newcommand{\intz}{\mathbb{Z}}

\newcommand{\ftl}{\tilde{f}}
\newcommand{\htn}{\hat{\theta}_{n}}

\makeatother

\begin{document}
\begin{frontmatter}

\title{Goodness of fit tests for a class of Markov random field models}
\runtitle{Goodness of fit tests for spatial data}

\begin{aug}
\author[A]{\fnms{Mark S.} \snm{Kaiser}\ead[label=e1]{mskaiser@iastate.edu}},
\author[B]{\fnms{Soumendra N.} \snm{Lahiri}\thanksref{t1}\ead[label=e2]{snlahiri@stat.tamu.edu}}
\and
\author[A]{\fnms{Daniel~J.}~\snm{Nordman}\corref{}\thanksref{t2}\ead[label=e3]{dnordman@iastate.edu}}
\runauthor{M. S. Kaiser, S. N. Lahiri and D. J. Nordman}
\affiliation{Iowa State University, Texas A\&M University and Iowa
State University}
\address[A]{M. S. Kaiser\\
D. J. Nordman\\
Department of Statistics\\
Iowa State University\\
Ames, Iowa 5001\\
USA\\
\printead{e1}\\
\phantom{E-mail: }\printead*{e3}}
\address[B]{S. N. Lahiri\\
Department of Statistics\\
Texas A\&M University\\
College Station, Texas 77843\\
USA\\
\printead{e2}} 
\end{aug}

\thankstext{t1}{Supported in part by NSF Grants DMS-07-07139
and DMS-10-07703.}

\thankstext{t2}{Supported in part by NSF Grant DMS-09-06588.}

\received{\smonth{4} \syear{2010}}
\revised{\smonth{11} \syear{2011}}

%
\begin{abstract}
This paper develops goodness of fit statistics that can be used to
formally assess Markov random field models for spatial data, when the
model distributions are discrete or continuous and potentially
parametric. Test statistics are formed from generalized spatial
residuals which are collected over groups of nonneighboring spatial
observations, called concliques. Under a hypothesized Markov model
structure, spatial residuals within each conclique are shown to be
independent and identically distributed as uniform variables. The
information from a series of concliques can be then pooled into
goodness of fit statistics. Under some conditions, large sample
distributions of these statistics are explicitly derived for testing
both simple and composite hypotheses, where the latter involves
additional parametric estimation steps. The distributional results are
verified through simulation, and a data example illustrates the method
for model assessment.
\end{abstract}

%

\end{frontmatter}

\section{Introduction}\label{sec1}

Conditionally specified models formulated on the basis of an
underlying Markov random field (MRF) are an attractive alternative
to continuous random field specification for the analysis of
problems that involve spatial dependence structures. By far
the most common of such models are those formulated using
a conditional Gaussian distribution (e.g., \cite{RH05}), but models may
also be constructed using a number of other conditional
distributions such as a beta \cite{HY08,KCL02},
binary \cite{CK09}, Poisson \cite{B74} or Winsorized Poisson
\cite{KC97}, and general specifications are
available for many exponential families \cite{ACS92,KCL02}.

In an applied spatial setting, we assume that observations are
available at a finite set of geo-referenced locations
$\{ \bfs_i\dvtx i=1, \ldots, N\}$, and to these locations we assign
the random variables $\{Y(\bfs_i)\dvtx i=1, \ldots, N\}$. In general,
locations are arbitrarily indexed in $d$-dimensional real space.
A MRF is typically constructed by specifying
for each location $\bfs_i$ a \textit{neighborhood},
consisting of other locations on which the full
conditional distribution of $Y(\bfs_i)$ will be functionally
dependent. Let the conditional cumulative distribution function (c.d.f.)
of $Y(\bfs_i)$ given $\{Y(\bfs_j)=y(\bfs_j)\dvtx j \neq i \}$ be denoted
as $F_i$ and define $\N_i \equiv\{\bfs_j \neq\bfs_i$,
and $F_i$ depends functionally on $y(\bfs_j) \}$.
Also define $\bfy(\N_i) \equiv\{ y(\bfs_j)\dvtx\bfs_j \in\N_i\}$.
The Markov assumption implies that
%
%
\begin{equation}
\label{eqn1}
F_i\bigl(\cdot| \{ y(\bfs_j)\dvtx\bfs_j \neq\bfs_i\}\bigr)
= F_i\bigl(\cdot|
\{y(\bfs_j)\dvtx\bfs_j \in\N_i\}\bigr) = F_i(\cdot| \bfy(\N_i)).
\end{equation}

A model is formulated by specifying,
for each $i=1, \ldots, N$, a conditional c.d.f.
in (\ref{eqn1}). Conditions necessary for a set of such
conditionals to correspond to a joint distribution for
$\{ Y(\bfs_1), \ldots, Y(\bfs_N) \}$ are given by
Arnold, Castillo and Sarabia \cite{ACS92} and a constructive process
with useful conditions sufficient for existence of a joint
are laid out in Kaiser and Cressie \cite{KC00}. Models may
be constructed for both discrete and continuous random
variables, on regular or irregular lattices, with or
without an equal number of neighbors for each location (including $\N_i
= \varnothing$ for some locations) and possibly including information
from spatial covariates.
The construction of models for applications is thus very flexible.

A number of our results and, in particular, Theorem \ref{Theorem2.1} to
follow, can
be generalized to some of the variable situations just described, but
it will be beneficial for developing theoretical results to define a
setting that is broad but highly structured. We desire a spatial
\textit{process} defined on grid nodes of the $d$-dimensional integer lattice
$\intz^d$, where $\intz= \{0,\pm1,\pm2,\ldots\}$. We stipulate a
number of restrictions for this process that, while not capable of
covering all of the finite-dimensional models mentioned previously, is
flexible enough to be meaningful in many applied situations.
We formally consider specifying an MRF model for a spatial process
$\bfY\equiv\{ Y(\bfs)\dvtx\bfs\in\intz^d\}$,
rather than a model~(\ref{eqn1}) developed with respect to a finite
collection of (possibly nonlattice)
data sites $\{Y(\bfs_i)\dvtx i=1, \ldots, N\}$. To this end, assume that
for any $\bfs\in\intz^d$ neighborhoods can be constructed using a
standard template $\M\subset\intz^d\setminus\{\bfzero\}$ as $\N(\bfs
) = \bfs+ \M$, with $|\M| < \infty$ denoting the size of~$\M$. Some
examples of $\M$ are given in the next section. We then assume that the
process $\bfY$ has a stationary distribution function $F(\cdot|\cdot)$
such that, for any \mbox{$\bfs\in\mathbb{Z}^d$}, the conditional c.d.f.
of $Y(\bfs)$ given all remaining variables $\{ Y(\bft)\dvtx\bft\in
\mathbb{Z}^d,\break \bft\neq\bfs\}$ can be written as
%
%
\begin{equation}
\label{processmodel}
F\bigl(\cdot| \{ Y(\bft)\dvtx\bft\in\mathbb{Z}^d, \bft\neq\bfs\}\bigr) =
F\bigl(\cdot| \{ Y(\bft)\dvtx\bft\in\N(\bfs) \}\bigr)
\end{equation}
under a Markov assumption.

Given a hypothesized or estimated model, our concern is how one might
conduct a goodness
of fit (GOF) procedure, either through informal diagnostics or by
using formal probability results that lead to a GOF test. The approach
we propose here may be viewed within
either the context of a~pure GOF test to address
the question of whether a (possibly fitted) model provides an adequate
description of observed data. This is an issue of
model assessment and different from model selection,
which has been considered, for example, with penalized pseudo-likelihood
for parametric MRF models; cf. \cite{CT06,GY99,JS96}. Additionally,
while other GOF tests may be possible
for certain joint model specifications (e.g., a frequency-domain approach
for Gaussian processes; cf. \cite{A93}), we focus solely on
conditional model specifications.
The GOF variates introduced
in the next section may be used as either diagnostic
quantities or as the basis for a formal GOF test as presented in
Section \ref{sec3}.

The remainder of this article is organized as follows. In
Section \ref{sec2} we introduce the concept of a conclique
and derive GOF variates that form
the basis of our approach, using an adaptation of a
multivariate probability integral transform (PIT). Section \ref{sec3}
develops a formal
methodology for combining these variates over concliques
to create GOF tests of Markov models under both
simple and composite hypotheses.
These tests are omnibus in the sense that they assess the hypothesized
model in total, including the neighborhood structure selected,
specification of dependence as isotropic or directional, and the form
of the modeled conditional distributions.
Theoretical results are presented in Section \ref{sec4} that
establish the limiting sampling distributions of GOF tests under the
null hypothesis.
Section \ref{sec5} describes a~numerical study
to support the theoretical findings. Section \ref{sec6} provides
an application of the GOF tests in model assessment for agricultural trials.
Section \ref{sec7} contains concluding remarks and discussions on extensions.
Section~\ref{sec8} provides a proof
of the foundational conclique result (Theorem \ref{Theorem2.1}), and
all other
proofs regarding the asymptotic
distribution of GOF test statistics appear in supplementary material
\cite{KLN11}.

\section{Generalized spatial residuals}\label{sec2}

In this section we derive the basic quantities that form the
basis for our GOF procedures. We consider these
quantities to be a type of generalized residuals because they fit
within the framework suggested by
Cox and Snell \cite{CS71}.
In particular, these generalized spatial residuals
will be derived using an extended version of Rosenblatt's \cite{R52}
multivariate PIT combined with a partitioning of
spatial locations into sets such that the residuals within each
set constitute a random sample from a uniform distribution
on the unit interval, under the true model. As discussed by
Brockwell \cite{Br07} and Czado et al. \cite{Cz09}, the
PIT formulation allows arbitrary model distributions to be considered
in assessing GOF, rather than simply
continuous ones. Similar transformations, with subsequent formal or
informal checks for uniformity,
have been important in evaluating the GOF of, and the quality of
predictive forecasts from,
various models for time series; cf. \cite{Da57,Da84,Di98,Gn07,JW07,Ju97}.

\subsection{Concliques}\label{sec2.1}

Before providing the transform that defines our generalized spatial
residuals, it is necessary to develop a method for partitioning
the total set of spatial locations at which observations are
available into subsets with certain properties. We call such
sets \textit{concliques} because they are defined as the converse of
what are called cliques by Hammersley and Clifford~\cite{HC71}. In the
case of regular lattices with neighborhoods defined using either
four-nearest or eight-nearest neighbor structures, concliques
correspond exactly to the so-called coding sets of Besag \cite{B74},
which were suggested for use in forming conditional likelihoods
for estimation. The key property of concliques, however,
allows construction of such sets in more general settings
including irregular lattices and hence the new name.

As defined in \cite{HC71}, a \textit{clique}
is a set of locations such that each location in the set
is a neighbor of every other location in the set.
Similar terminology exists in graph theory, where a subset
of graph vertices (e.g., locations) form a clique if
every two vertices in the subset are connected by an edge \cite{S78}.
We define a
\textit{conclique} as a set of locations such that no location
in the set is a~neighbor of any other location in the set.
Any two members of a conclique may share common neighbors,
but they cannot be neighbors themselves.
Additionally, every set of a single location can be treated as both a
clique or conclique.
In the parlance of graphs, the analog of a conclique
is a so-called ``independent set,'' defined by a set of vertices in
which no
two vertices share an edge. This particular graph terminology
conflicts with the probabilistic notion of independence, while a
``conclique''
truly represents
a \textit{conditionally} independent set of locations
in a MRF model.

While the result of the next subsection holds for any collection
of concliques, in practice what is desired is a collection of
concliques that suitably partition all observed locations.
To achieve this under the process model~(\ref{processmodel}), we
identify a collection of concliques
$\{\C_j\dvtx j=1, \ldots, q\}$ that partition the entire grid $\intz^d$.
We define a collection of concliques to be a~\textit{minimal conclique
cover} if it contains the smallest number of concliques needed to
partition the set of all locations. In graph theory, this concept is related
to determining the smallest (or chromatic) number of colors needed to
color a graph (with no two edge-connected vertices sharing the same color)
or, equivalently, the smallest number of independent sets needed to
partition graph vertices \cite{W96}.
In practice, identifying a minimal conclique cover is valuable since
our procedure produces one test statistic for each conclique in a
collection, and those statistics must then be combined into one overall
value for a formal GOF test.
%
%
\begin{Example}[(A $4$-nearest neighbor model on $\intz^2$)]\label{Example2.1}
Here, let $\bfs=(u,v)' \in\intz^2$ for a horizontal coordinate $u$
and a vertical coordinate $v$.
The neighborhood structure of a $4$-nearest neighbor model is produced
with the template $\M=\{(-1,0)', (1,0)',(0,1)',
(0,-1)'\}$, so that $\N(\bfs)$ for a given location $\bfs$ and
neighbors $\ast$ is as shown
in the following figure:
\[
\matrix{
\cdot&\ast&\cdot\cr
\ast& \bfs&\ast\cr
\cdot&\ast&\cdot}
\]
In this case, the minimal conclique cover contains two members, $\C_1$
and $\C_2$,
with elements denoted by $1$'s and $2$'s, respectively, as shown below.

\begin{center}
\textit{Minimal conclique collection for a $4$-nearest neighbor model}:
\end{center}
\[
\matrix{
1&2&1&2&1&2&1&2&1&2\cr
2&1&2&1&2&1&2&1&2&1\cr
1&2&1&2&1&2&1&2&1&2\cr
2&1&2&1&2&1&2&1&2&1\cr
1&2&1&2&1&2&1&2&1&2}
%
\]
\end{Example}
%
%
\begin{Example}[(An $8$-nearest neighbor model on $\intz^2$)]\label{Example2.2}
As in the previous example, let $\bfs= (u,v)'$ but take $\M= \{
(u,v)'\dvtx\max\{ |u|, |v| \} = 1 \}$. The neighborhood structure of
an $8$-nearest neighbor model is then shown in the following figure for
a location $\bfs\in\intz^2$ and neighbors $*$:
\[
\matrix{
\ast&\ast&\ast\cr
\ast& \bfs&\ast\cr
\ast&\ast&\ast}
\]
For the $8$-nearest neighbor model, there are four concliques in the
minimal cover,
$\C_1,\ldots,\C_4$, with elements denoted by $1$'s, $2$'s, $3$'s and $4$'s
in the following figure, respectively.

\begin{center}
\textit{Minimal conclique cover for an $8$-nearest neighbor model}:
\end{center}
\[
\matrix{
1&2&1&2&1&2&1&2&1&2\cr
3&4&3&4&3&4&3&4&3&4\cr
1&2&1&2&1&2&1&2&1&2\cr
3&4&3&4&3&4&3&4&3&4\cr
1&2&1&2&1&2&1&2&1&2}
%
\]
\end{Example}

\subsection{Defining generalized spatial residuals}\label{sec2.2}

$\!\!\!$Let $\{A(\bfs) \dvtx\bfs\in\mathbb{Z}^d\}$ denote a~collection of
independent and identically distributed (i.i.d.)
random variables, which
are Uniform $(0,1)$ and also independent of the spatial process~$\bfY$.
For any $\bfs\in\intz^d$, we then define a random generalized
spatial residual as
%
%
\begin{eqnarray}
\label{genresid}
U(\bfs) & =& \bigl(1-A(\bfs)\bigr) \cdot F\bigl(Y(\bfs)|\{Y(\bft)\dvtx\bft
\in\N(\bfs) \} \bigr)\nonumber\\[-8pt]\\[-8pt]
&&{} + A(\bfs) \cdot F^{-}\bigl(Y(\bfs)|\{Y(\bft)\dvtx\bft
\in\N(\bfs) \} \bigr),\nonumber
\end{eqnarray}
where $F(\cdot|\cdot)$ denotes the (stationary) c.d.f. from
(\ref{processmodel}), and $F^{-}(\cdot|\cdot)$ denotes the left limit
of the c.d.f., that is,
$F^{-}(y | \{ Y(\bft)\dvtx\bft\in\N(\bfs) \}) = P( Y(\bfs) < y |\break \{
Y(\bft)\dvtx\bft\in\N(\bfs) \})$, $y\in\mathbb{R}$.
This residual applies the notion of a randomized PIT \cite{Br07},
allowing for a noncontinuous c.d.f.
$F(\cdot|\cdot)$ to be considered. When $F(\cdot|\cdot)$ is
continuous, the spatial residual
reduces to a PIT $U(\bfs) = F(Y(\bfs)|\{Y(\bft)\dvtx\bft\in\N(\bfs)
\} )$ in Rosenblatt's \cite{R52} format.
Given that a collection of concliques is available for a
particular situation, the fundamental result that serves as
the basis for our GOF procedures is as follows.
%
%
\begin{Theorem}\label{Theorem2.1}
Let the spatial process $\{ Y(\bfs)\dvtx\bfs\in
\mathbb{Z}^d \}$ have conditional distribution functions as in
(\ref{processmodel}), and let $\{ \mathcal{C}_j\dvtx j=1, \ldots, q\}$
be a
collection of concliques that partition the integer grid $\intz^d$.
Then for any $j=1, \ldots, q$, the variables $\{U (\bfs)\dvtx\bfs\in
\mathcal{C}_j\}$ given by (\ref{genresid}) are i.i.d.
Uniform $(0,1)$ variables.
\end{Theorem}

Typically, the conditional c.d.f. $F(\cdot|\cdot)$
of expression (\ref{processmodel}) will be a parameterized function,
and we
now write this as $F_{\te}( \cdot| \cdot)$ to emphasize
the parametrization. Let $\te_0$ denote the true value
of the parameter. In an application we have available
a set of observations
taken to represent realizations of the random variables
$\{Y(\bfs_i)\dvtx i=1, \ldots, N\}$. Theorem \ref{Theorem2.1} indicates that
if we compute generalized spatial residuals as, in the notation of
(\ref{genresid}),
%
%
\begin{eqnarray}
\label{realizedgenresid}
U(\bfs_i) &=& \bigl(1-A(\bfs_i)\bigr)\cdot F_{\theta_0}\bigl(y(\bfs_i) | \{
y(\bft)\dvtx
\bft\in\N(\bfs_i) \} \bigr) \nonumber\\[-8pt]\\[-8pt]
&&{} + A(\bfs_i) \cdot F^{-}_{\theta_0}\bigl(y(\bfs_i)| \{
y(\bft)\dvtx\bft\in\N(\bfs_i) \} \bigr),\qquad
\bfs_i \in\mathcal{C}_j,\nonumber
\end{eqnarray}
then within any conclique $\mathcal{C}_j$ these variables
should behave as a random sample from a uniform distribution
on the unit interval. If we use a minimal conclique cover
having $q$ members, then we will have $q$ sets of residuals,
each of which should behave as a random sample from a uniform
distribution. These sets of residuals will not, however, be
independent, so we will not have a total collection that behaves
as $q$ independent random samples.

In practice we will usually also replace the parameter $\theta$
with an estimate~$\hat{\theta}$ computed
on the basis of the observations so that, technically,
the values within any conclique will not actually be independent
either. We expect, however, that if the model is appropriate, then
these residuals will exhibit approximately the same behavior
as independent uniform variates,
in the same way that ordinary residuals from a linear regression
model with normal errors behave as an approximate random
sample of normal
variates, despite the fact that they cannot technically
represent such a sample.

A basic diagnostic plot can be constructed by plotting the
empirical distribution function of each set of residuals
$\{u(\bfs_i)\dvtx\bfs_i \in\C_j \}$, $j=1, \ldots, q$,
and examining them for departures from
a standard uniform distribution function. See, for instance,
Gneiting et al. \cite{Gn07}, Section 3.1, for a summary
of graphical approaches for exploring uniformity in PIT values. Tests
for uniformity
may be used for individual sets of residuals to guide the
decision about whether a given fitted model is adequate or to
choose between two competing (even nonnested) models. Such procedures
do not constitute a~formal
GOF test, however, because there is no guarantee
that results will agree across differing sets of residuals in a
conclique cover.
Formal procedures for combining evidence from the
residual sets into one overall GOF test are
presented in the next section.

\section{Methodology: Goodness of fit tests}\label{sec3}

\subsection{General setting}\label{sec3.1}

Suppose that for a set of locations on the $d$-dimen\-sional integer
lattice $\{\bfs_1, \ldots, \bfs_N\} \subset\intz^d$,
we want to assess the GOF
of a~conditional model specification, based on a set of observed values
$\{Y(\bfs_i)\dvtx i=1, \ldots, N \}$. We assume that the observed values
are a partial realization of a~class of process models defined on $\intz
^d$ for which the conditional
c.d.f. of~$Y(\bfs)$ given $\{Y(\bft)\dvtx\bft\neq\bfs\}$ belongs to a
class of parameterized conditional distribution functions,
%
%
\begin{equation}
\label{hypothclass}
{\mathcal F}_{\theta} =
\bigl\{ F_{\theta}\bigl( \cdot| \{ Y(\bft)\dvtx\bft\in\N(\bfs) \}\bigr)\dvtx
\theta\in\Theta\bigr\},
\end{equation}
where
$\Theta\subseteq\rl^p$, $1\leq p<\infty$, is a parameter space, $\N(\bfs
) = \bfs+ \M$ and, analogously to (\ref{processmodel}), $\M\subset
\intz^d \setminus\{ \bfzero\} $. Two testing problems fit into this
framework, where the null hypothesis is simple and where it is composite.

In the next subsections, we describe GOF tests for simple and composite
hypotheses based on the observations $\{ Y(\bfs_i)\dvtx i=1, \ldots, N\}$,
which are assumed to have arisen in the following way.
Suppose that $R \subset\rl^d$ denotes a~sampling region within which
$N$ observations are obtained at a set of sampling locations $\cS_N
\equiv R \cap\intz^d = \{ \bfs_1, \ldots, \bfs_N\}$. Define the
interior of the set of sampling locations~as $\cS^{\mathrm{int}}_N \equiv\{
\bfs\in\cS_N\dvtx\N(\bfs) \subset\cS_N \}$. Locations in this set
are those sampling locations for which all neighbors are also sampling
locations, allowing generalized spatial residuals to be computed for
all $\bfs\in\cS^{\mathrm{int}}_N$, even if the physical sampling region $R$ is
irregular.
Finally, let $\C_{1N},\ldots,\C_{qN}$ denote the conclique partition of
$\cS^{\mathrm{int}}_N$ determined by $\C_{jN} = \C_j \cap\cS_N^{\mathrm{int}}$, $j=1,
\ldots, q$. In practice we will desire a minimal conclique cover but
this is not necessary in what follows.

\subsection{Testing a simple null hypothesis}\label{sec3.2}

First consider the case of the simple $(S)$ null hypothesis in which
the testing problem is
given by, for some specified $\theta_0 \in\Theta$,
\begin{eqnarray*}
&&H_0(S)  \mbox{: The data } \{Y(\bfs_i)\dvtx i=1, \ldots, N\} \mbox{
represent a partial sample of } \\
&&\hphantom{H_0(S)\mbox{: }} \mbox{the process model class (\ref{hypothclass}) with } \theta
=\theta_0; \\
&&H_1(S) \mbox{: Not $H_0(S)$}.
\end{eqnarray*}

To construct test statistics appropriate for these hypotheses, we
consider the generalized spatial residuals under $H_0(S)$,
%
%
\begin{eqnarray}
\label{gresidsimple}
U(\bfs) &=& \bigl(1-A(\bfs)\bigr) \cdot F_{\theta_0} \bigl( Y(\bfs) | \{ Y(\bft
)\dvtx\bft\in\N(\bfs) \}\bigr)\nonumber\\[-8pt]\\[-8pt]
&&{} + A(\bfs) \cdot F^{-}_{\theta_0} \bigl( Y(\bfs) | \{
Y(\bft)\dvtx\bft\in\N(\bfs) \}\bigr),\qquad \bfs\in\cS_N^{\mathrm{int}}.\nonumber
\end{eqnarray}
Now define, for $j=1,\ldots,q$, the (generalized
residual) empirical distribution function over the $j$th conclique by
\[
G_{jN}(u)=\frac{1}{|\C_{jN}|}\sum_{\bfs\in\C_{jN}}\ind\bigl(U(\bfs)\leq u\bigr),
\]
$u\in[0,1]$. Here and in the following, $\ind(A)$ denotes the indicator
function of a~statement $A$, where $\ind(A) = 1$ if $A$ is true and
$\ind(A) = 0$ otherwise.
Note that under $H_0(S)$, $E \{ G_{jN}(u) \} =u$,
$u\in[0,1]$, as a result of Theorem \ref{Theorem2.1}. Hence, to assess
the GOF of
the model over the $j$th conclique $\C_j$, we consider the scaled deviations
of the empirical distribution function from the Uniform $(0,1)$ distribution,
%
%
\begin{equation}
\label{wjnsimple}
W_{jN}(u)\equiv N^{1/2}\bigl(G_{jN}(u) - u\bigr),\qquad u\in[0,1].
\end{equation}

A number of GOF test statistics for testing $H_0(S)$
may be obtained by combining the $W_{jN}$'s
in different ways:
%
%
\begin{eqnarray}
\label{teststat1simple}
T_{1N} &=& {\max_{j=1,\ldots,q} \sup_{u\in[0,1]}} |W_{jN}(u)|,\\
\label{teststat2simple}
T_{2N}&=& \Biggl( \frac{1}{q}\sum_{j=1}^q\Bigl[\sup_{u\in[0,1]}
|W_{jN}(u)|
\Bigr]^2\Biggr)^{1/2},\\
\label{teststat3simple}
T_{3N} &=& \max_{j=1,\ldots,q} \biggl(\int_0^1|W_{jN}(u)|^r
\,du\biggr)^{{1/r}}, \\
\label{teststat4simple}
T_{4N}&=& \frac{1}{q}\sum_{j=1}^q
\biggl(\int_0^1|W_{jN}(u)|^r\,du\biggr)^{{1/r}},
\end{eqnarray}
where $r\in[1,\infty)$ in (\ref{teststat3simple}) and (\ref{teststat4simple}).
Note that $T_{1N}$ and $T_{2N}$ are obtained by combining conclique-wise
Kolmogorov--Smirnov test statistics, while $T_{3N}$ and $T_{4N}$ are
obtained by
combining conclique-wise Cram\'er--von Mises test statistics.
While our statistics are based exclusively on paired differences
(e.g., $G_{jN}(u)-u$, $u\in[0,1]$), other test statistics may be
formulated to assess agreement
between the empirical $G_{jN}$ and Uniform$(0,1)$ distributions,
such as GOF tests based on $\phi$-divergences studied in \cite{JW07}.
In Section \ref{sec4}, we provide asymptotic distributions for the empirical
processes (\ref{wjnsimple}),
which may be an ingredient for determining limit distributions of
statistics based on $\phi$-divergences; cf. Theorem 3.1 \cite{JW07}.

\subsection{Testing a composite null hypothesis}\label{sec3.3}

The composite ($C$) null hypothesis can be stated as
\begin{eqnarray*}
&&H_0(C)\mbox{: The data } \{Y(\bfs_i)\dvtx i=1, \ldots, N\} \mbox{
represent a partial sample of } \\
&&\hphantom{H_0(C)\mbox{: }} \mbox{some member of the process model class (\ref{hypothclass})
for an unknown $\theta$}; \\
&& H_1(C) \mbox{: Not $H_0(C)$}.
\end{eqnarray*}

Let $\hth$ denote an estimator of $\te$ based on
$\{Y(\bfs_i)\dvtx i=1, \ldots, N \}$. Since $\te$ is unknown,
instead of the $U(\bfs)$'s of (\ref{gresidsimple}), we
work with an estimated version of the generalized spatial residuals,
%
%
\begin{eqnarray}
\label{gresidcomposite}
\hU(\bfs)&=& \bigl(1-A(\bfs)\bigr) \cdot F_{\hth} \bigl( Y(\bfs) | \{ Y(\bft)\dvtx
\bft\in\N(\bfs) \}\bigr)\nonumber\\[-8pt]\\[-8pt]
&&{} + A(\bfs) \cdot F^{-}_{\hth} \bigl( Y(\bfs) | \{ Y(\bft
)\dvtx\bft\in\N(\bfs) \}\bigr),\qquad \bfs\in\cS_N^{\mathrm{int}},\nonumber
\end{eqnarray}
where, as before, $\N(\bfs)= \bfs+ \M$.
Note that if $\hth$ is a
reasonable estimator of $\te$ and if
$F_{\te}(\cdot|\cdot)$ is a smooth function of $\te$,
then the $\hU(\bfs)$'s of (\ref{gresidcomposite}) are approximately
distributed as
Uniform $(0,1)$. This suggests that we can base
tests of $H_0(C)$ versus $H_1(C)$ on the processes
%
%
\begin{equation}
\label{wjncomposite}
\hW_{jN}(u)\equiv N^{1/2}\bigl(\hG_{jN}(u) - u\bigr),\qquad u\in[0,1],
\end{equation}
for $j=1,\ldots,q$, where
\[
\hG_{jN}(u) =\frac{1}{|\C_{jN}|}\sum_{\bfs\in\C_{jN}}\ind\bigl(
\hU(\bfs) \leq u\bigr),\qquad u\in[0,1].
\]
The test statistics for testing $H_0(C)$ versus $H_1(C)$ are
now given by
%
%
\begin{equation}
\label{eqn310}
\hT_{1N},\ldots, \hT_{4N},
\end{equation}
where $\hT_{jN}$ is obtained by replacing $W_{jN}$ in expressions
(\ref{teststat1simple})--(\ref{teststat4simple}) with
$\hW_{jN}$.
In the next section, we describe the limit distributions of the
test statistics under the null hypothesis.

\section{Asymptotic distributional results}\label{sec4}

\subsection{Basic concliques}\label{sec4.1}

To formulate large sample distributional results for the
GOF statistics, we shall assume that the concliques $\C_1,\ldots,\C_q$
used for these statistics can be ``built up'' from unions of
structurally more basic concliques, say $\C_1^*,\ldots,\C_{q^*}^*$,
$q^* \geq q$.
For any given template $\M\subset\intz^d\setminus\{\mathbf{0}\}$
defining neighborhoods as $\N(\bfs) = \bfs+ \M$, $\bfs\in\intz^d$,
we suppose such concliques are constructed as follows.

Let
$\mathbf{e}_i\in\mathbb{Z}^d$ denote a vector with $1$ in the $i$th component
and 0 elsewhere, and define $m_i \equiv\max\{ |\mathbf{e}_i^\prime
\mathbf{s}|\dvtx\mathbf{s}\in\mathcal{M} \}$ as the maximal absolute value
of $i$th component\vadjust{\goodbreak} over integer vectors in the neighborhood template
$\mathbf{s}\in\mathcal{M}$, $i=1,\ldots,d$.
Define a collection of sublattices as
%
%
\begin{equation}\label{eqn-ccl-str}
\C_j^* =\{\bfa_j + \Delta\bfs\dvtx\bfs\in\intz^d\},\qquad
j=1,\ldots, q^* \equiv\prod_{i=1}^d (m_i+1),
\end{equation}
where $\Delta= \operatorname{diag}(m_1 +1,\ldots, m_d+1)$ is a positive
diagonal matrix and
\[
\mathbf{a}_j\in\mathcal{I}\equiv\{(a_1,\ldots,a_d)^\prime\in\mathbb
{Z}^d \dvtx0 \leq a_i \leq m_i, i=1,\ldots,d\},
\]
where $\mathbf{a}_j \neq\mathbf{a}_k$ if $\C^*_j \neq\C^*_k$.

Proposition \ref{Proposition4.1} shows that these sets provide a
collection
of ``basic'' concliques (or coding sets) since locations within the
same sublattice $\C_j^*$ are separated by directional distances $\Delta$
that prohibit neighbors within $\C_j^*$.
Additionally,\vspace*{1pt} the proposition gives a simple rule for
merging basic concliques~$\C_j^*$
to create larger concliques~$\C_j$.
In the following, write $\pm\mathcal{M}=\mathcal{M}\cup-\mathcal{M}$,
and define $\|\mathbf{s}\|_\infty= \max_{1 \leq i \leq d}|s_i|$ for
$\mathbf{s}=(s_1,\ldots,s_d)^\prime\in\mathbb{Z}^d$.
%
%
\begin{Proposition}\label{Proposition4.1}
Under the process assumptions of
Theorem \ref{Theorem2.1} and for any neighborhood specified by a finite subset
$\M\subset\intz^d\setminus\{\bfzero\}$:

\begin{longlist}[(a)]
\item[(a)]
sets $\C_1^*,\ldots, \C^*_{q^*}$ of form (\ref{eqn-ccl-str})
are concliques that partition $\intz^d$;

\item[(b)] if $\mathbf{a}_1,\ldots, \mathbf{a}_{i},\mathbf{a}_{i+1}\in\mathcal
{I}$, $i \geq1$, such that $\mathcal{C}\equiv\bigcup_{j=1}^i \mathcal
{C}^*_{j}$ is a conclique, then $\mathcal{C} \cup\C^*_{i+1}$ is a
conclique if and only if
\[
\mathbf{a}_j-\mathbf{a}_{i+1} + \Delta\mathbf{s} \notin\pm\mathcal
{M} \qquad\mbox{for all $\mathbf{s}\in\mathbb{Z}^d$, $\|\mathbf{s}\|
_\infty\leq1$, and any $j=1,\ldots,i$}.
\]
\end{longlist}
\end{Proposition}

In addition to providing a systematic approach for building concliques,
the purpose of this basic conclique representation is to allow the covariance
structure of the limiting Gaussian process of the conclique-wise
empirical processes [cf. (\ref{wjnsimple})] to be written explicitly
and to simplify the distributional results to follow
(as basic concliques $\C_j^*$ above have a uniform structure and are
translates of one another). With many
Markov models on a regular lattice described by the neighborhoods in
Besag \cite{B74} involving coding sets
or ``unilateral'' structures, there is typically no loss of generality
in building a~collection of concliques $\C_1, \ldots, \C_q$ from such
basic concliques.
We illustrate Proposition \ref{Proposition4.1} with some examples.
{\renewcommand{\theExample}{2.\arabic{Example}}
\setcounter{Example}{0}
%
\begin{Example}[(Continued)]
Under the four-nearest neighbor
structure in $\mathbb{Z}^2$, we have $\mathcal{M}= \{ \pm(0,1)^\prime,
\pm(1,0)^\prime\} = \pm\mathcal{M}$, $m_1=m_2=1$, $\Delta=
\operatorname{diag}(2,2)$ and \mbox{$q^*=4$}, so there are four basic
concliques $\{C_j^*\}
_{j=1}^4$ determined by the vectors
\[
\mathbf{a}_1=(0,0)^\prime,\qquad \mathbf{a}_2=(1,1)^\prime,\qquad
\mathbf{a}_3=(1,0)^\prime, \qquad\mathbf{a}_4=(0,1)^\prime.
\]
Because $\mathbf{a}_2 - \mathbf{a}_1 + \Delta\cdot\mathbf{s} =
(1,1)^\prime+ 2 \mathbf{s} \notin\pm\mathcal{M}$
for any $\mathbf{s}\in\mathbb{Z}^2$, $\|\mathbf{s}\|_\infty\leq1$,
then $\mathcal{C}_{1}\equiv\mathcal{C}_{1}^*\cup\mathcal{C}_{2}^*$ is
a conclique,\vadjust{\goodbreak} and, similarly, so is $\mathcal{C}_2\equiv\mathcal
{C}_{3}^*\cup\mathcal{C}_{4}^*$.
Additionally, Proposition~\ref{Proposition4.1} shows also that $\mathcal
{C}_1$, $\mathcal{C}_2$
cannot be further merged so that these represent the previously
illustrated minimal conclique cover.
\end{Example}
%
%
\begin{Example}[(Continued)]
Under the eight-nearest
neighbor structure in $\mathbb{Z}^2$, we have that $\mathcal{M}= \{ \pm
(0,1)^\prime, \pm(1,0)^\prime, \pm(1,1)^\prime, \pm(1,-1)^\prime\}$
and the basic concliques $\{\C_j^*\}_{j=1}^4$ are the\vspace*{1pt} same as
in Example \ref{Example2.1} and correspond to Besag's~\cite{B74} coding sets.
However, these basic concliques cannot be merged into larger concliques by
Proposition \ref{Proposition4.1} and hence match the minimal cover of
four concliques as
illustrated previously (i.e., $\C_j=\C^*_j$).
\end{Example}}
%
\setcounter{Example}{0}
\begin{Example}\label{Example2.3}
Under a ``simple unilateral'' neighbor
$\mathcal{M}= \{ (0,-1)^\prime$, $(-1,0)^\prime\}$ in $\mathbb{Z}^2$
(cf. \cite{B74}, Section 6.2), the basic concliques are again the same
and
Proposition \ref{Proposition4.1} gives $\mathcal{C}_{1}\equiv\mathcal
{C}_{1}^*\cup
\mathcal{C}_{2}^*$, $\mathcal{C}_2\equiv\mathcal{C}_{3}^*\cup\mathcal
{C}_{4}^*$ as a minimal conclique cover.
\end{Example}

\subsection{Asymptotic framework}\label{sec4.2}

We now consider a sequence of sampling regions $R_n$ indexed by $n$.
For studying the large sample properties of the proposed
GOF statistics, we adopt an ``increasing domain
spatial asymptotic'' structure \cite{C93}, where the
sampling region $R_n$ becomes unbounded
as $\nti$.
Let
$R_0$ be an open connected subset of $(-1/2,1/2]^d$ containing
the origin.
We regard $R_0$
as a ``prototype'' of the sampling region $R_n$. Let
$\{\lambda_n\}_{n\geq1}$ be a
sequence of positive numbers such that $\lambda_n
\raw\infty$
as $n\rightarrow\infty$.
We assume that the sampling region $R_n=\lambda_nR_0$ is
obtained by ``inflating'' the set $R_0$ by the scaling
factor $\lambda_n$ (cf. \cite{LKCH99}).
Since the origin is assumed to lie in $R_0$, the shape of $R_n$ remains
the same for different values of $n$.
To avoid pathological cases, we assume that for any sequence of real
numbers $\{a_n\}_{n\geq1}$ with $a_n\rightarrow0+$ as $n\rightarrow
\infty$, the
number of cubes of the lattice $a_n\intz^d$ that intersect both
$R_0$ and $R_0^c$ is $O((a_n)^{-(d-1)})$ as $n\rightarrow\infty$.
This implies that, as the sampling region grows, the number of
observations near the boundary of $R_n$ is of smaller order $O(N_n^{(d-1)/d})$
than the total number $N_n$ of observations in~$R_n$ so that the volume
of $R_n$, $N_n$ and the number of interior locations are equivalent as
$n \rightarrow\infty$.
The boundary condition on~$R_0$ holds for
most regions~$R_n$ of practical
interest, including common convex subsets of $\rl^d$,
such as rectangles and ellipsoids, as well as for many
nonconvex star-shaped sets in~$\rl^d$. (Recall that a set
$A\subset\rl^d$ is called star-shaped if for any $x\in A$, the
line segment joining $x$ to the origin lies in $A$.)
The latter class of sets may have a
fairly irregular shape. See, for example,
\cite{L99,SC94} for more details.

We want to assess the GOF
of the process model specification (\ref{processmodel}),
under either the simple or composite hypothesis sets of Section \ref{sec3}.
As described in Section~\ref{sec3.1}, we suppose that the spatial
process is
observed at
locations on the integer grid $\intz^d$ that fall
in the sampling region $R_n$ producing a set of sampling locations $\cS
_{N_n}$ (indexed by $n$).
To simplify notation,
we will use $\cS_n$ rather than the more cumbersome $\cS_{N_n}$ and
$\cS_n^{\mathrm{int}}$ rather than $\cS_{N_n}^{\mathrm{int}}$.
Similarly, we\vadjust{\goodbreak} will use $W_{jn}$ to denote the empirical distribution of
generalized spatial residuals for the $j$th conclique under a simple
hypothesis as given by (\ref{wjnsimple}) with $N=N_n$ and $T_{1n},
\ldots, T_{4n}$, the corresponding test statistics of
(\ref{teststat1simple})--(\ref{teststat4simple}). Also, $\hat{W}_{jn}$,
and $\hat{T}_{1n}, \ldots, \hat{T}_{4n}$ will denote the quantities in
(\ref{wjncomposite}) and~(\ref{eqn310}) with $N=N_n$.

\subsection{Results for the simple testing problem}\label{sec4.3}

For studying the asymptotic distribution of the test statistics
$T_{1n},\ldots,T_{4n}$ under the null hypothesis~$H_0(S)$, we shall make
use of the following condition, which imposes the structure on the
concliques described in Section \ref{sec4.1}.
\begin{longlist}[\textit{Condition} (C.1):]
\item[\textit{Condition} (C.1):] Each conclique $\C_1,\ldots,\C_q$ is union of basic concliques
$\C_1^*,\ldots,$ $C_{q^*}^*$ as in (\ref{eqn-ccl-str}). Namely, for each
$j=1,\ldots,q$,
there exists $\mathcal{J}_j \subset\{1,\ldots,q^*\equiv
\operatorname{det}(\Delta)\}$ where $\C_j= \bigcup_{i \in\mathcal{J}_j}
\C_i^*$ and
the index sets $\{\mathcal{J}_j\}_{j=1}^q$ are disjoint.
\end{longlist}

The following result gives the asymptotic null distribution of
conclique-wise empirical processes $\bfW_n =(W_{1n},\ldots, W_{qn})'$
based on the scaled and centered empirical distributions $W_{jn}(u)$,
$u\in[0,1]$, as in (\ref{wjnsimple}). Note that, while each individual
empirical process $W_{jn}$ can be expected to weakly converge to a
Brownian bridge under $H_0(S)$ (cf. \cite{B03}), the limit law of
$\bfW_n$ will not similarly be distribution-free due to the dependence
in observations across concliques. In particular, the null model
$F_{\theta_0}$ influences the asymptotic covariance structure of
$\bfW_n$.\vspace*{1pt}

Let ${\mathcal L}_{\infty}^q$ denote the collection of bounded
vector-valued functions $\bff=(f_1,\ldots,$ $f_q)'\dvtx[0, 1] \rightarrow
\rl
^q$ defined on the unit interval. Also, let $|B|$ denote the size of a
finite set $B \subset\mathbb{R}$.
%
%
\begin{Theorem}\label{Theorem4.2}
Suppose that condition \textup{(C.1)} holds.
Then, there exists a
zero-mean vector-Gaussian process
$\bfW(u) = (W_{1}(u),\ldots, W_{q}(u))', u\in[0,1]$,
with continuous sample paths on $[0,1]$ (with probability
1)
such that
\[
\bfW_n \stackrel{d}{\raw} \bfW\qquad\mbox{as $\nti$}
\]
as elements of ${\mathcal L}_{\infty}^q$.
Further, $P(\bfW(u)=\bfzero)=1$ for $u=0,1$ and
the $q\times q$ covariance matrix function of $\bfW$ is given by
\[
E W_j(u)W_k(v) = \cases{
\displaystyle \frac{\operatorname{det}(\Delta)}{|\mathcal{J}_j|}(\min\{u, v\} - uv ),
&\quad if $j=k$,\vspace*{2pt}\cr
\displaystyle \frac{\operatorname{det}(\Delta)}{|\mathcal{J}_j| \cdot
|\mathcal{J}_k|}\sum_{i \in\mathcal{J}_j, l \in\mathcal{J}_k } \sigma
_{i,l}(u,v), &\quad if $j\neq k$,}
\]
for $0\leq u,v\leq1$, $1\leq j,k\leq q$ and
\begin{eqnarray*}
\sigma_{i,l}(u,v) &\equiv&
\sum_{\mathbf{s}\in\mathbb{Z}^d, \|\mathbf{s}\|_\infty\leq1} \{
P[ U(\mathbf{0})\leq u, U(\mathbf{a}_l-\mathbf{a}_i + \Delta\mathbf
{s})\leq v ]-uv\}\\
&&\hspace*{48.5pt}{}\times\mathbb{I}(\mathbf{a}_l-\mathbf{a}_i +
\Delta\mathbf{s} \in\pm\mathcal{M}).
\end{eqnarray*}
\end{Theorem}

The indicator function $\mathbb{I}(\cdot) $ above pinpoints terms in
the covariance expression which
automatically vanish by the independence of residual
variables~$U(\mathbf{s})$ within conclique structures (Theorem \ref{Theorem2.1}).
For example, when it is possible to combine two basic concliques $\C
_i^*$ and $\C_l^*$, $i \neq l$, into a~larger conclique, Proposition
\ref{Proposition4.1} gives that, for all $\|\mathbf{s}\|_\infty\leq1$,
it holds that
$\mathbf{a}_i-\mathbf{a}_l + \Delta\mathbf{s}
\notin\mathcal{M}$ and so above $\mathbb{I}(\mathbf{a}_i-\mathbf{a}_l
+ \Delta\mathbf{s} \in\pm\mathcal{M})=0$.
All sums in the limiting covariance structure then involve only a
finite number of terms.

As a direct implication of Theorem \ref{Theorem4.2}, we get the
following result
on the asymptotic null distribution of the test statistics
$T_{1n},\ldots,
T_{4n}$.
%
%
\begin{Corollary}\label{Corollary4.3}
Under the conditions of Theorem \ref{Theorem4.2},
\[
T_{jn}\stackrel{d}{\raw} \varphi_j(\bfW) \qquad\mbox{as $\nti$}
\]
for $j=1,\ldots,4$, where
the functionals' $\varphi_j$'s are defined by
%
%
\begin{eqnarray}\label{eqn-functionals}
\varphi_{1}(\bff) &=& {\max_{1\leq j \leq q} \sup_{u\in
[0,1]}}|f_{j}(u)|,\nonumber\\
\varphi_{2}(\bff) &=& \biggl(\frac{1}{q}
\sum_{1\leq j \leq q} \Bigl[{\sup_{u\in[0,1]}}
|f_{j}(u)|\Bigr]^2\biggr)^{1/2}, \nonumber\\[-8pt]\\[-8pt]
\varphi_{3}(\bff) &=& \max_{1\leq j \leq q} \biggl(\int_0^1|f_{j}(u)|^r
\,du\biggr)^{1/r},\nonumber\\
\varphi_{4}(\bff) &=& \frac{1}{q}\sum_{1\leq j \leq q}
\biggl(\int_0^1|f_{j}(u)|^r \,du \biggr)^{1/r}\nonumber
\end{eqnarray}
for $\bff=(f_1,\ldots,f_q)' \in{\mathcal L}_{\infty}^q$,
and for a given $r\in[1,\infty)$.
\end{Corollary}

\subsection{Results for the composite testing problem}\label{sec4.4}

As for the simple testing problem, here we first derive the asymptotic
null distribution
of the conclique-wise empirical processes
$\hat{\mathbf{W}}_n =(\hat{W}_{1n},\ldots, \hat{W}_{qn})'$
based on scaled and centered empirical distributions $\hat
{W}_{jn}(u)$, $u\in[0,1]$, in (\ref{wjncomposite}).

Note that the estimator $\htn$
appears in each summand in $\hW_{jn}$ through the estimated generalized
spatial residuals (\ref{gresidcomposite}). In such situations,
a common standard approach to deriving asymptotic distributions of
empirical processes is based on the concept of \textit{uniform asymptotic
linearity}
in some local neighborhood of the true parameter value $\te_0$
(cf. \cite{K70,VW96}). However, this
approach is \textit{not} directly applicable here due to the form of the
conditional distribution functions in
(\ref{hypothclass}) when considered as functions of $\theta\in\Theta
$. To establish the limit distribution, we embed the
empirical process of the estimated generalized residuals
in an enlarged space, namely,
the space of locally bounded $q$-dimensional vector functions on $[0,1]$,
equipped with the metric of uniform convergence on compacts,
and then
use a version of the continuous mapping theorem; the argument details
are provided in~\cite{KLN11}.

We require some notation and conditions in addition to those introduced
in the earlier section.
Letting again
$|B|$ denote the size of a finite set $B$, define the strong mixing
coefficient of the process $\{Y(\bfs)\dvtx\bfs\in\intz^d \}$
by
\begin{eqnarray*}
\alpha(a;b) &=& \sup\{|P(A\cap B) - P(A)P(B)|\dvtx A\in\D(S_1), B\in\D(S_2),
\\
&&\hspace*{43pt}
|S_1|\leq b,
|S_2|\leq b, d(S_1,S_2) \geq a, S_1, S_2 \subset\intz^d\},
\end{eqnarray*}
where $\D(S)=\sigma\langle Y(\bfs)\dvtx\bfs\in S \rangle$
generically denotes the $\sigma$-algebra generated by variables $Y(\bfs
)$ with locations in $S \subset\intz^d$,
$d(S_1,S_2)=\inf\{\|\bfs-\bft\|_1\dvtx\bfs\in S_1,\break\bft\in S_2\}$, $\|
\mathbf{x}\|_1=\sum_{i=1}^d |x_i|$ for $\mathbf{x}=(x_1,\ldots
,x_d)^\prime\in\mathbb{R}^d$, and $P(\cdot)$ represents probabilities
for the process.
Write
$F^{(1)}_\te(\cdot| \cdot)$ and $F^{(1)-}_\te(\cdot| \cdot)$ to
denote $p\times1$ vectors of first order partial derivatives
of $F_\te(\cdot| \cdot)$ and $F_\te^{-}(\cdot| \cdot)$ with respect
to $\te$, when these exist. Let $U_\theta(\bfzero) = (1-A(\bfzero
))\cdot F_{\te}(Y(\bfzero)| \{Y(\bft)\dvtx\bft\in\M\}) +A(\bfzero)\cdot
F^{-}_{\te}(Y(\bfzero)| \{Y(\bft)\dvtx\bft\in\M\}) $, and
denote $\bfU_\theta^{(1)}(\bfzero)\in\mathbb{R}^p$ as the vector of
partial derivatives
of $U_\theta(\bfzero)$ with respect to $\te$, when this exists.

\textit{Condition} (C.2):
\begin{longlist}
\item There exist constants $\delta_0 \in(0,1)$, $c_0\in
(0,\infty)$ such that
\[
\bigl| P\bigl( U_{\te}(\bfzero) \leq u \bigr) - P\bigl(
U_{\te_0}(\bfzero) \leq v \bigr) \bigr| \leq c_0 [ \|\theta-\theta_0 \|
+|u-v|]
\]
for all $0 \leq u,v \leq1$ and $\theta\in\Theta$ satisfying $\max\{\|
\theta-\theta_0 \|, |u-v|\}\leq\delta_0$.\vspace*{1pt}

\item $\sup\{ \|F^{(1)}_\te( y | \bfx)\| + \|F^{(1)-}_\te( y | \bfx)\|
\dvtx
\|\theta-\theta_0 \| \leq\delta_0, y\in\mathbb{R},\bfx\in\mathbb
{R}^p \}
\leq c_0$.

\item $ E \{ \sup_{\|\te-\te_0\|<\delta} \|\bfU_\theta^{(1)}(\bfzero)
- \bfU_{\te_0}^{(1)}(\bfzero) \| \} = o(\delta)$ as
$\delta\rightarrow0$.
\end{longlist}

\textit{Condition} (C.3):
Suppose that the joint distribution of
$(U_{\te_0}(\bfzero), \bfU_{\te_0}^{(1)}(\bfzero))$ is absolutely continuous
with respect to
$L\times\mu$ with
Radon--Nikodym derivative $\ftl(u,\bfx)$, where $L$ is the Lebesgue
measure on $\rl$, and
$\mu$ is a $\sigma$-finite measure on $\rl^p$. Suppose that
\[
\lim_{t\raw\infty} \sup_{u\in(0,1)}\int_{\|\bfx\|>t}
\|\bfx\| \ftl(u,\bfx)\,d\mu(\bfx) =0
\]
and
\[
\int\|\bfx\| \cdot\sup\{|\ftl(u,\bfx) - \ftl(v,\bfx)|\dvtx|u-v|\leq
\delta\}
\,d\mu(\bfx) \raw0
\]
as $\delta\raw0+$.

\textit{Condition} (C.4):
\begin{longlist}
\item There exist zero-mean random variables
$\{\bfV(\bfs)\dvtx
\bfs\in\intz^d\}$
such that
\[
N_n^{1/2}(\htn-\te_0) = N_n^{-1/2}\sum_{\bfs\in\cS_n}\bfV(\bfs)
+o_p(1).
\]
\item For each $\bfs\in\intz^d$, the variable
$\bfV(\bfs)=(V_1(\bfs),\ldots,V_p(\bfs))^\prime$ is $\D(\bfs+\M
)$-measurable.\vadjust{\goodbreak}

\item
There exist $a \in(2,\infty)$, $\kappa>0$ such that
$\sup\{E\|\bfV(\bfs)\|^{2+\kappa}\dvtx\bfs\in\intz^d\}<\infty$ and
\[
\sum_{j=1}^\infty j^{d-1}\alpha(j
;1)^{{\kappa}/({2+\kappa})}<\infty,\qquad
\sum_{j=1}^\infty j^{d(2r-1)}\alpha(j ;2r-1)^{{1}/{a}}<\infty
\]
for some integer $r$ satisfying $r>(p+1)/(1-a^{-1}) $.

\item $\Sigma\equiv\lim_{\nti} \operatorname{Var}( N_n^{-1/2}\sum_{\bfs\in\cS_n}\bfV
(\bfs))$
exists and is nonsingular.
\end{longlist}

Conditions (C.2) and (C.3) are exclusively used for handling the
effects of the perturbation of the empirical process of
the generalized residuals due the estimation of $\te$.
The first displayed condition in (C.3)
is an uniform integrability condition, while the second one
is a continuity condition on the densities $\ftl(\cdot,\cdot)$
(in~$u$) in a weighted $L^1(\mu)$-norm. Without loss of generality, we
shall suppose
that $\ftl(u,\bfx)=0$ for all $u \notin(0,1)$ except on a set of $\bfx
$-values with $\mu$-measure zero.
Condition (C.4)
allows us to relate the limit law of the (unperturbed)
empirical process part with the variability in estimating
$\te$ by $\htn$.
If the conditional
model specification is such that the spatial process satisfies Dobrushin's
uniqueness condition (cf. \cite{G95}), then the MRF is
strongly mixing (actually, $\phi$-mixing) at an exponential
rate and, hence, mixing conditions in (C.4)
trivially hold.
%
%
\begin{Theorem}\label{Theorem4.4}
Suppose that conditions
\textup{(C.1)--(C.4)} and the composite null hypothesis $H_0(C)$
hold. Then, there exist a zero-mean vector-Gaussian process
$\bfW(u) = (W_{1}(u),\ldots, W_{q}(u))', u\in[0,1]$,
with continuous sample paths on $[0,1]$ (with probability
1) and a random variable $\bfZ=(Z_1,\ldots,Z_p)^\prime\sim N_p(\bfzero
,\Sigma)$, both defined
on a common probability space,
such that as $n \rightarrow\infty$,
\[
\hat{\mathbf{W}}_n \stackrel{d}{\raw} \bfW+ \bfone\cdot\bfZ'\int
\bfx\ftl(\cdot,\bfx)\,d\mu(\bfx)
\]
as elements of ${\mathcal L}_{\infty}^q$,
where $\bfone=(1,\ldots,1)'\in\rl^q$.
The $q\times q$ covariance matrix function of $\bfW$ is
as in Theorem \ref{Theorem4.2} and for $j=1,\ldots,q$, $k=1,\ldots,p$
and $u\in(0,1)$,
\[
E W_j(u)Z_k= \frac{1}{|\mathcal{J}_j|} \sum_{i \in\mathcal{J}_j} \sum
_{\bfs\in\intz^d}
E\bigl(V_k(\bfs- \bfa_i)\cdot
\ind\bigl(U(\bfzero)\leq u\bigr)\bigr).
\]
\end{Theorem}

The following result is a direct consequence of
Theorem \ref{Theorem4.4} and gives the asymptotic distribution of the
test statistics
under the composite null~$H_0(C)$.

%
\begin{Corollary}\label{Corollary4.5}
Under the conditions of Theorem \ref{Theorem4.4},
\[
\hT_{jn} \stackrel{d}{\raw}
\varphi_j \biggl(\bfW+ \bfone\cdot\bfZ'\int\bfx\ftl(\cdot,\bfx)\,d\mu
(\bfx)\biggr)\qquad
\mbox{as $\nti$}
\]
for $j=1,\ldots,4$, where the functionals $\varphi_{j}$'s are as
defined in\vadjust{\goodbreak}
(\ref{eqn-functionals}).
\end{Corollary}

Under the composite null $H_0(C)$, the limiting distributions involved
are not distribution-free (i.e., depending on the true model c.d.f.
$F_{\te_0}$ in a complex covariance structure). Empirical processes
based on PIT residuals with parameter estimates are known to exhibit
this behavior in other inference scenarios with time series and
independent data (cf. \cite{D75}), and often two general approaches are
considered for implementing GOF tests \cite{K05}: resampling or
Khmaladze's \cite{K81} martingale transformation. The latter involves
a~type of continuous de-trending to minimize effects of parameter
estimation and has been applied to obtain asymptotically
distribution-free tests with other model checks using residual
empirical processes based on estimated parameters (cf.
\cite{K93,KK04}). In particular, Bai \cite{B03} justified this
transformation for tests in parametric, conditionally specified
(continuous) distributions for time series, but considered only one
empirical process of residuals. If modified to the spatial setting,
this result would entail a transformation of~$\hW_{jN}$ from one
conclique $j=1,\ldots,1$ so that its limiting distribution is Brownian
motion and distribution-free under~$H_0(C)$. The complication here is
that with residual empirical processes from multiple concliques, after
applying a conclique-wise transformation, the resulting limit
distribution of a~test statistic under~$H_0(C)$ would not be
distribution-free due to dependence across concliques (akin to Theorem
\ref{Theorem4.2} in the case of no parameter estimation). Another option might be to
use plug-in estimates of the covariance structure, using, for example,
that asymptotic variances of maximum likelihood and pseudolikelihood
estimators (i.e., $\Sigma$ in Theorem \ref{Theorem4.4}) are known for
some Markov field models \cite{GY99}. But one would also have to
estimate other complicated covariances in the limiting distribution of
Theorem \ref{Theorem4.4}, which might be possible with subsampling
variance estimation \cite{SC94}.

Spatial resampling methodologies,
such as the block bootstrap (cf. \cite{L03}, Chapter~12), might also
be used
to approximate sampling distributions of GOF statistics based on
spatial residuals and knowledge of
the limit distributions in Theorem \ref{Theorem4.4} could be applied to toward
justifying such bootstrap estimators.
Simulations in Section \ref{sec5} also suggest that
the finite sample versions of the GOF statistics appear to converge
fairly quickly to their limits, at least in the case of simple null
hypotheses. This implies that, in application, large-sample bootstrap
approximations of finite-sample sampling distributions
may be reasonable.
The theoretical development of a spatial bootstrap for our GOF
statistics is outside of the scope of this paper,
but in Section \ref{sec6} we use a parametric spatial bootstrap to calibrate
GOF test statistics for a composite null hypothesis.

\section{Numerical results}\label{sec5}
Here we provide a small numerical verification of the large sample
distributional results in the simple null hypothesis case, considering
observations generated from a conditional Gaussian MRF on $\mathbb
{Z}^2$ with a four-nearest neighbor structure
specified by $\mathcal{M}=\{\pm(0,1)', \pm(1,0)'\}$ as in Example
\ref{Example2.1}.
The conditional model family (\ref{hypothclass}) of
$Y(\mathbf{s})$ given $\{Y(\bft)\dvtx\break\bft\in\N(\bfs) \}$, $\mathbf{s}\in
\mathbb{Z}^2$ ($\mathcal{N}(\mathbf{s})=\mathbf{s}+\mathcal{M}$), is
normal with mean
$\mu_{\alpha,\eta}(\mathbf{s}) \equiv\alpha+ \eta\sum_{\bft\in
\mathcal{N}(\mathbf{s})}[Y(\bft)-\alpha]$,
and variance $\tau^2>0$, where $E(Y(\mathbf{s}))=\alpha\in\mathbb{R}$
is the marginal process mean and $|\eta| < 0.25$ denotes a dependence
parameter. In total, the model parameters $\theta$ are
$(\alpha,\tau,\eta)^\prime$.

\subsection{Limit distributions under a simple null hypothesis}\label{sec5.1}

We first examine asymptotic null distributions of GOF test statistics
in the simple testing problem
$H_0(S)\dvtx(\alpha, \tau, \eta)^\prime= (\alpha_0,\tau_0,\eta_0)^\prime$
of Section \ref{sec3.2} with residuals (\ref{gresidsimple}) given by
$U(\mathbf
{s}) = \Phi[\{Y(\mathbf{s})-\mu_{\alpha_0,\eta_0}(\mathbf{s}) \}/\tau
_0] $. Here $\Phi(\cdot)$ denotes the standard normal cumulative
distribution function, and, for simplicity, we will write hypothesized
parameters $\alpha_0,\tau_0,\eta_0$ as $\alpha,\tau,\eta$ in the following.

As described in Section \ref{sec3.1}, the four-nearest-neighbor structure
produces a minimal cover of two concliques
$\C_1, \C_2$ (cf. Example \ref{Example2.1}), each of which is a union
of two basic
concliques $\C^*_1,\ldots,\C^*_4$
provided in Section \ref{sec4.1}.
These concliques yield an empirical distribution process $\bfW_{n} =
(W_{1n}, W_{2n})^\prime$
and GOF test statistics $T_{1n},\ldots,T_{4n}$ as in
(\ref{teststat1simple})--(\ref{teststat4simple}). By Theorem
\ref{Theorem4.2}, $\bfW_{n}$
has a mean-zero Gaussian limit $\bfW=(W_1,W_2)^\prime$ with covariances
%
%
\begin{equation}
\label{eqnsim}
E W_{j}(u) W_k(v) = \cases{
2 (\min\{u,v\}-uv), &\quad if $j=k$,\cr
8\bigl[ P\bigl(X_1 \leq\Phi^{-1}(u), X_2 \leq\Phi^{-1}(v) \bigr) - uv\bigr], &\quad if
$j \neq k$,}\hspace*{-35pt}
\end{equation}
$u,v\in[0,1], j,k\in\{1,2\}$, where vectors $(X_1,X_2)$ in
(\ref{eqnsim}) are bivariate normal, with marginally standard normal
distributions and correlation $-\eta$. Hence, under the simple null
hypothesis, the limit process depends on $(\alpha, \tau, \eta)^\prime$
only through the dependence parameter $\eta$, which we denote by writing
$\bfW\equiv\bfW_\eta$.

%
\begin{figure}

\includegraphics{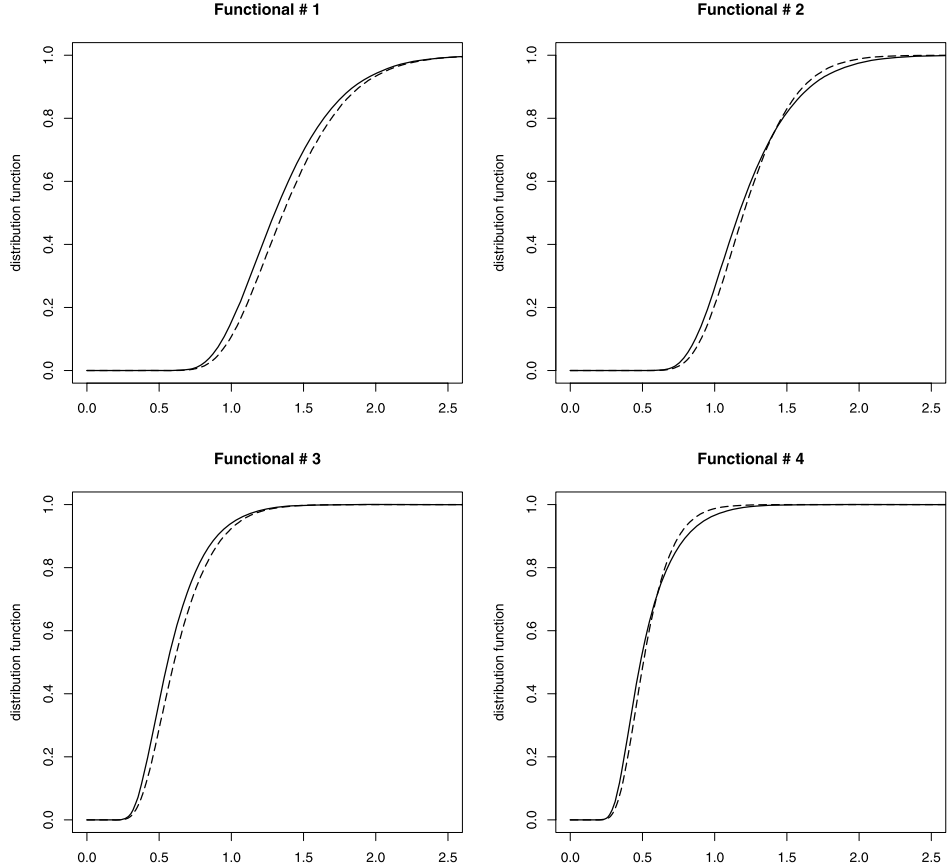}

\caption{Cumulative distribution functions
$F_{\varphi_j(\mathbf{W}_\eta)}(w)= P[ \varphi_j(\mathbf{W}_\eta)\leq w
]$, $w \in\mathbb{R}$, for limit functionals
$\varphi_1(\mathbf{W}_\eta),\ldots,\varphi_4(\mathbf{W}_\eta)$ for
$\eta=0.1$ (dashed) and $\eta=0.24$ (solid).}
\label{fig1}
\end{figure}

To understand the distribution of $\varphi_j(\bfW_\eta)$, $j=1,2,3,4$,
as the asymptotic limit of GOF statistics $T_{jn}$ under Corollary
\ref{Corollary4.3},
we simulated from the theoretical Gaussian process $\bfW_\eta$ as follows.
For each value of $\eta=0,0.1,0.24$, we generated 50,000 sequences of
mean-zero bivariate Gaussian variables $(W_1(i/3001)$, $W_2(i/3001))$,
$i=0,\ldots,3001$, with covariance structure (\ref{eqnsim}) over a
grid in $[0,1]$; the sequence length of $3002$ was dictated by
computational stability.
These provide approximate observations of $\bfW_\eta$, with
$\eta$ values chosen to reflect no, weak and strong forms of positive
spatial dependence.
Cumulative distribution functions of each functional $\varphi_1(\mathbf
{W}_\eta),\ldots,\varphi_4(\mathbf{W}_\eta)$ were then approximated from
$\mathbf{W}_\eta$-realizations. The resulting distribution curves
appear in Figure \ref{fig1} for $\eta=0.1$ and $\eta=0.24$, with
$\varphi
_3(\mathbf{W}_\eta)$ and $\varphi_4(\mathbf{W}_\eta)$ computed using
$r=2$ in (\ref{eqn-functionals}).

\subsection{Comparisons to finite sample distributions}\label{sec5.2}

\begin{figure}

\includegraphics{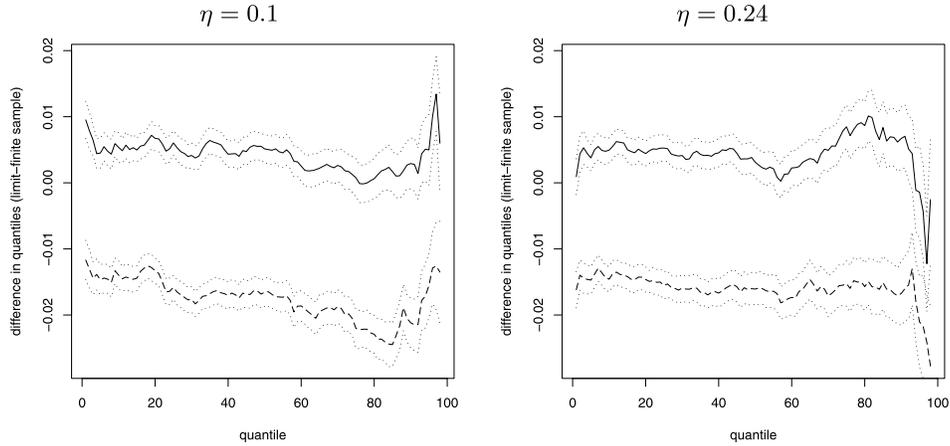}

\caption{Difference in quantiles for $\varphi_2(\bfW_\eta)$ and
$T_{2N}$ when $N=100$ (dashed line) and 900 (solid line) for
$\eta=0.1,0.24$. Pointwise 95\% confidence bands (dotted) indicate the
Monte Carlo error in each difference.}
\label{fig2}
\end{figure}

To compare the agreement of finite sample distributions of $T_{jN}$
under the simple null hypothesis with their limit distributions $\varphi
_j(\mathbf{W}_\eta)$, $j=1,\ldots,4$, we simulated samples on two grid
sizes, a $10 \times10$ grid having $N=100$ locations and a $30 \times
30$ grid having $N=900$.
Here, we simulated 50,000 realizations of conditional Gaussian samples
(setting $\alpha=0$ and $\tau=1$ with no loss of generality) and
evaluated functionals $T_{1N},\ldots,T_{4N}$ to approximate the
finite-sample distributions of these GOF statistics. Figure \ref{fig2}
shows the
difference between
the quantiles of the limit $\varphi_2(\mathbf{W}_\eta)$ and those of
$T_{2N}$ for $\eta=0.1$ and $\eta=0.24$;
the agreement among quantiles for functional 2 is quite good even
though this plot
was one exhibiting the largest quantile-mismatches among the four GOF
functionals.
Table \ref{table1} shows the proportion of GOF statistics $T_{jN}$
falling above
the 95th and 99th quantiles
of the corresponding limit $\varphi_j(\bfW_\eta)$ distribution,
$j=1,2,3,4$. The agreement between the finite-sample
and theoretical limit distributions is again close in Table \ref{table1}.

%

%
\begin{table}[b]
\caption{Proportion of GOF statistics $T_{jN}$ from a
conditional Gaussian model falling above the 95th and 99th
quantiles (denoted $q_{95,\eta}$ and $q_{99,\eta}$) of the their limit
$\varphi_j(\bfW_\eta)$ distribution, $j=1,2,3,4$, for sample sizes
$N=100$ and $N=900$ and with
dependence parameters $\eta=0,0.1,0.24$}\label{table1}
\begin{tabular*}{\tablewidth}{@{\extracolsep{\fill}}lccccccccc@{}}
\hline
&\multicolumn{1}{c}{}& \multicolumn{4}{c}{\textbf{\% of} $\bolds{T_{jN} > q_{95,\eta}}$}
& \multicolumn{4}{c}{\textbf{\% of} $\bolds{T_{jN} > q_{99,\eta}}$} \\ [-4pt]
&& \multicolumn{4}{c}{\hrulefill} & \multicolumn{4}{c@{}}{\hrulefill}\\
$\bolds{\eta}$ & $\bolds{N}$ & $\bolds{j=1}$ & $\bolds{2}$ & $\bolds{3}$
& $\bolds{4}$ & $\bolds{j=1}$ & $\bolds{2}$ & $\bolds{3}$ & $\bolds{4}$ \\
\hline
0 &100 &4.67&4.38&5.09&4.97&0.90&0.90&0.98&0.91\\
0 &900 & 5.11&4.91&4.95&4.86&1.03&1.09&1.03& 1.03\\
0.1&100 & 4.60&4.60&4.88&4.92&0.94&0.95&0.95&1.08\\
0.1&900 & 5.11&5.13&5.05&5.08&1.08&1.13&1.07&1.15\\
0.24&100 & 4.52&4.57&5.02&5.06&0.80&0.76&0.86&0.92\\
0.24&900 & 4.92&4.97&4.86&5.03&0.97&0.96&0.93&0.95 \\
\hline
\end{tabular*}
\end{table}

For various sample sizes and dependence parameters, Table \ref{table2} compares
the finite-sample distributions of the four GOF statistics $\{T_{jN}\}
_{j=1}^4$ against their limiting distributions
$\varphi_j(\mathbf{W})$ in terms of a Kolmogorov--Smirnov $D_{\mathrm{KS}}$ and
a Cram\'er--von Mises-like $D_{\mathrm{CM}}$ distance metric, defined by
\begin{eqnarray*}
D_{\mathrm{KS}}(X,Z) &\equiv& \sup_{t\in\mathbb{R}}|F_X(t) - F_Z(t) |,\\
D_{\mathrm{CM}}(X,Z) &\equiv& \biggl[\int|F_X(t) - F_Z(t) |^2 \,dt\biggr]^{1/2},
\end{eqnarray*}
relative to the cumulative distributions $F_X,F_Z$ of arbitrary random
variables $X,Z$. To interpret the relative values of these metrics
in assessing~the distributional distance between $T_{jN}$ and $\varphi
_j(\bfW_\eta)$,
it is helpful to reference $D_{\mathrm{KS}},D_{\mathrm{CM}}$ values for comparing the
distributions of
$\varphi_j(\bfW_{\eta_1})$ and~$\varphi_j(\bfW_{\eta_2})$ over
parameters $\eta_1\neq\eta_2$, which Table \ref{table2} also provides.
%
%
\begin{table}
\caption{Computed values ($\times$1000) from distance
metrics comparing finite-sample distributions
of statistics $T_{1N},\ldots,T_{4N}$ to their limiting distributions
$\varphi_1(\mathbf{W}_{\eta}),\ldots,\varphi_4(\bfW_\eta)$}
\label{table2}
\begin{tabular*}{\tablewidth}{@{\extracolsep{\fill}}lcd{2.1}d{2.1}d{3.1}d{2.1}
d{2.1}d{2.1}d{2.1}d{2.1}@{}}
\hline
&& \multicolumn{4}{c}{$ \bolds{D_{\mathrm{KS}}(\varphi_j(\mathbf{W}_{\eta}),
T_{jN})}$} & \multicolumn{4}{c@{}}{$\bolds{ D_{\mathrm{CM}}(\varphi_j(\mathbf
{W}_{\eta}), T_{jN})}$} \\[-4pt]
&& \multicolumn{4}{c}{\hrulefill} &\multicolumn{4}{c@{}}{\hrulefill}\\
$\bolds{\eta}$ & $\bolds{N}$
& \multicolumn{1}{c}{$\bolds{j=1}$}
& \multicolumn{1}{c}{$\bolds{2}$} & \multicolumn{1}{c}{$\bolds{3}$}
& \multicolumn{1}{c}{$\bolds{4}$} & \multicolumn{1}{c}{$\bolds{j=1}$}
& \multicolumn{1}{c}{$\bolds{2}$}
& \multicolumn{1}{c}{$\bolds{3}$} & \multicolumn{1}{c@{}}{$\bolds{4}$}\\
\hline
0 &100 & 19.6 &23.0 & 8.9 & 9.1 & 12.9 &14.8 &3.7&3.2\\
0 &900 & 6.7 &9.7 &3.4 &4.8& 3.8 &4.7& 1.0 &1.4\\
0.1&100 & 24.0 &27.7& 5.3 & 5.4 & 16.5& 18.2& 2.2&1.7\\
0.1&900 & 9.9 &10.0 & 6.5 &5.7 & 5.3 &4.7& 2.5 &1.8\\
0.24&100 & 21.6 &25.2& 7.2 &7.0 & 14.7 &15.6 & 2.7 &2.2\\
0.24&900 & 9.1 &8.2& 4.2& 3.8 & 4.8 &4.8 &1.4&1.4\\
\hline
& &\multicolumn{4}{c}{$\bolds{D_{\mathrm{KS}}(\varphi_j(\mathbf{W}_{\eta_1}),
\varphi_j(\mathbf{W}_{\eta_2}))}$} & \multicolumn{4}{c@{}}{$\bolds{
D_{\mathrm{CM}}(\varphi_j(\mathbf{W}_{\eta_1}), \varphi_j(\mathbf{W}_{\eta
_2}))}$}\\[-4pt]
&& \multicolumn{4}{c}{\hrulefill} &\multicolumn{4}{c@{}}{\hrulefill}\\
$\bolds{\eta_1}$ & $\bolds{\eta_2}$ & \multicolumn{1}{c}{$\bolds{j=1}$}
& \multicolumn{1}{c}{\textbf{2}}
& \multicolumn{1}{c}{\textbf{3}} & \multicolumn{1}{c}{\textbf{4}}
& \multicolumn{1}{c}{$\bolds{j=1}$} & \multicolumn{1}{c}{\textbf{2}}
& \multicolumn{1}{c}{\textbf{3}} & \multicolumn{1}{c@{}}{\textbf{4}} \\
\hline
0 & 0.1\hphantom{0} & 14.4 &11.9& 14.9 & 13.4 & 7.0 & 5.7& 7.7& 6.1\\
0.1& 0.24 & 70.0& 59.6 &90.3 &72.9 & 49.9 &35.0 &50.4 &33.1\\
0& 0.24 &81.8 &69.1 &102.3 &84.3 & 56.6 & 40.5 & 58.0 & 38.9\\
\hline
\end{tabular*}
\end{table}
Generally, the convergence of the finite-sample distributions $T_{jN}$
to their limits $\varphi_j(\mathbf{W}_\eta)$ appears to occur fairly
uniformly over different dependence parameters $\eta$ and, relative to
the distributional differences
among different limits [e.g.,~$\varphi_j(\bfW_{\eta_1})$ and~$\varphi
_j(\bfW_{\eta_2})$], the agreement in distributions of $T_{jN}$ and
$\varphi_j(\mathbf{W}_\eta)$
is quite close even for samples of size $100$.

\subsection{Power of GOF statistics under simple null hypothesis}\label{sec5.3}

\begin{figure}

\includegraphics{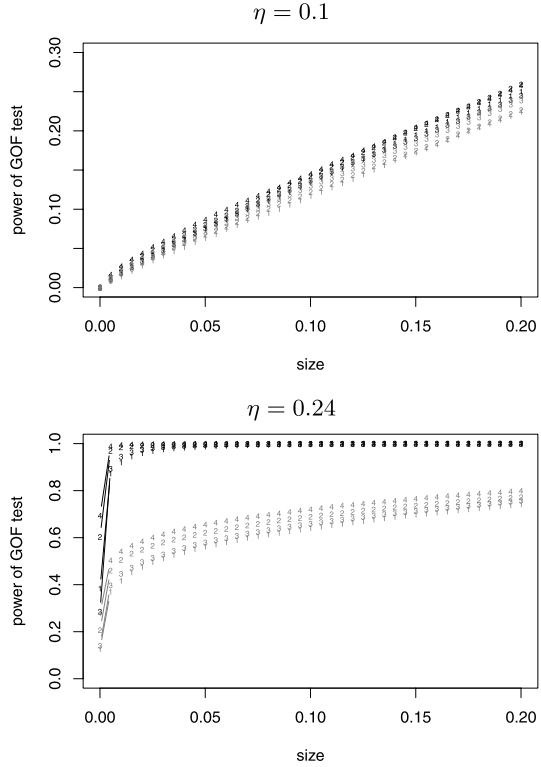}

\caption{Plots of power versus size $\gamma$ for GOF tests of
$H_0\dvtx\eta=0$ in conditional Gaussian models (fixed $\alpha=0,\tau=1$)
based on functionals $T_{1N},\ldots,T_{4N}$, determined by data
generated under $\eta=0.1,0.24$. In these power versus size curves,
each functional is numbered 1--4 under sample sizes $N=100$ (grey) and
$N=900$ (black).}\label{fig3}
\end{figure}

Under the simple null $H_0(S)\dvtx(\alpha,\tau,\eta)^\prime
=(0,1,0)^\prime
$, we next consider the power of GOF tests based on
statistics $T_{1N},\ldots,T_{4N}$ computed from conditional Gaussian
data generated with $\eta=0.1$ and $\eta=0.24$ and $\alpha=0,\tau=1$.
This gives an idea of the power in testing a hypothesis of no spatial
dependence, when the data exhibit
forms of positive dependence, both fairly weak ($\eta=0.1$) and strong
($\eta=0.24$).
For a given GOF statistic $T_{jN}$ from a sample of size $N=100$ or
$N=900$, a size $\gamma$ test is conducted by rejecting $H_0$ if
$T_{jN}$ exceeds the $1-\gamma$ quantile of the limit distribution
$\varphi(\mathbf{W}_{\eta=0})$ under the null hypothesis.
Figure \ref{fig3} plots power versus size $\gamma$ for these tests when
$\eta
=0.1$
%
%
and $\eta=0.24$, based on 50,000 simulated data sets. Power is
low under the alternative $\eta=0.1$, as might be expected, but
considerably higher when $\eta=0.24$.
Tests with functionals $T_{2N}, T_{4N}$ (based conclique-wise averages
of GOF statistics) tend to perform similarly and exhibit slightly
more power than tests with~$T_{1N}, T_{3N}$ (based conclique-wise
maxima of GOF statistics).

\section{An application to agricultural field trials}\label{sec6}
\subsection{The problem}\label{sec6.1}

Besag and Higdon \cite{BH99} present an analysis of six agricultural
field trials of corn varieties conducted in North Carolina using
a~hierarchical model that included an intrinsic Gaussian MRF as an
improper prior for spatial structure. An intrinsic Gaussian MRF results
from fixing dependence parameters at the boundary of the parameter
space. In discussion of this paper, Smith \cite{S99} raised the
question of what diagnostics were available to examine potential
evidence for spatial structure based on the available data, and
presented variograms of three of the trials. Kaiser and Caragea \cite
{KC09} used data from these same three trials to illustrate a
model-based diagnostic they called the $S$-value. Questions about the
spatial structure suggested by the data included the possibilities of
nonstationarity and directional dependencies. Here, we use data from
all six trials to examine the question of whether a simple model with
constant mean and unidirectional dependence can be rejected as a
plausible representation of spatial structure. Our question is simply
one of whether a basic Gaussian MRF with constant mean and a single
dependence parameter could be rejected as a possible data generating
mechanism for the data, not whether it might be be most preferred model
available.

Each field trial consisted of observations of yield from $64$ corn
varieties with each variety replicated $3$ times in each trial. The
spatial layout of each trial was essentially that of a $11 \times18$
regular lattice, although the last column of that lattice contained
only $5$ locations. After subtracting variety by trial means in the
same manner as \cite{KC09,S99}, we deleted the last column to obtain a
rectangular $11 \times17$ lattice containing $187$ observations for
each trial. We assumed a four-nearest-neighborhood structure but
without use of a border strip, so that locations had a variable number
of neighboring observations, $4$ for each of the $135$ interior
locations, $3$ for each of the $48$ edge locations, and $2$ for each of
the four corner locations.

\subsection{The model}\label{sec6.2}

Although each trial should nominally have marginal mean zero, to
examine a full composite setting we fit a model with conditional
Gaussian distributions having expected values $\{ \mu(\bfs_i)\dvtx i=1,
\ldots, n\}$ and constant conditional variance $\tau^2$ where, with
$N_i$ denoting the neighborhood of location $\bfs_i; i=1, \ldots, n$,
%
%
\begin{equation}
\label{eqnmsk}
\mu(\bfs_i) = \alpha+ \eta\sum_{\bfs_j \in N_i} \{ y(\bfs_j) - \alpha
\}.
\end{equation}
The joint distribution of this model is then Gaussian with marginal
means $\mbalpha$ an $n$-vector with each element equal to $\alpha$ and
covariance matrix $(I-C)^{-1} M$ where $I$ is the $n \times n$ identity
matrix, $M$ is an $n \times n$ diagonal matrix with all nonzero
entries equal to $\tau^2$ and $C=\eta H$ with $H$ an $n \times n$
matrix having element $(i,j)$ equal to $1$ if locations $\bfs_i$ and
$\bfs_j$ are neighbors and $0$ otherwise. With this structure, the
parameter space of $\eta$ can be determined to be $(-0.2563,
0.2563)$ based on eigenvalues of $H$ (cf. \cite{C93}); this differs
slightly from the parameter space for a lattice with
four-nearest-neighborhood structure wrapped on a torus due to the size
of the lattice and the use of varying numbers of neighbors for edge locations.

\subsection{The GOF procedure}\label{sec6.3}

The model of expression (\ref{eqnmsk}) was fit to (centered) data from
each of the six trials using maximum likelihood estimation. Generalized
spatial residuals were computed for each of the two concliques, one
having $93$ and the other $94$ locations. Using the fitted models, a
parametric bootstrap procedure was used to arrive at $p$-values for
each of the four test statistics introduced as $\hat{T}_{jN}; j=1,
\ldots, 4$, in Section \ref{sec3.3}. For each fitted model (i.e., trial) $5000$
bootstrap data sets were simulated using a Gibbs algorithm with a
burn-in of $500$ and spacing of $10$, which appeared adequate to result
in convergence of the chain based on scale reduction factors~\cite{GR}
and eliminate dependence between successive data sets based on
autocorrelations. Model (\ref{eqnmsk}) was fit, generalized spatial
residuals produced and the four test statistics computed for each
bootstrap data set, from which $p$-values were taken as the proportion
of simulated test statistic values greater than those from the actual
data sets. Bootstrap data sets were also used to produce percentile
bootstrap intervals for parameters (cf. \cite{DH97}). Percentile
intervals were chosen because basic bootstrap intervals extended beyond
the parameter space for $\eta$ for each of the six trials.

\subsection{Results}\label{sec6.4}

Results of estimation are presented in Table \ref{table3}. Intervals
were computed
at the $95\%$ level and values for $\eta$ are reported to four decimal
places because estimates tended to be close to the upper boundary of
the parameter space ($0.2563$). Overall, estimation was fairly similar
for these six trials, which were conducted in different counties of
North Carolina, including an indication of high variability in
estimating these parameters, particularly $\alpha$ and $\tau^2$.
Estimates of $\eta$ indicate moderate to strong spatial structure in
each of the six trials, and estimates of $\tau^2$ indicate substantial
local variability despite this structure.

%
\begin{table}
\caption{Estimates for conditional Gaussian models fit to data from six
agricultural field trials; the point estimates
for $\alpha$ for all trials differ from zero by at most $10^{-15}$}
\label{table3}
\begin{tabular*}{\tablewidth}{@{\extracolsep{\fill}}ld{3.2}cccc@{}}
\hline
& \multicolumn{2}{c}{\textbf{Point}}
& \multicolumn{3}{c@{}}{\textbf{Interval}} \\[-4pt]
& \multicolumn{2}{c}{\hrulefill}
& \multicolumn{3}{c@{}}{\hrulefill}\\
\textbf{Trial} & \multicolumn{1}{c}{$\bolds{\tau^2}$}
& \multicolumn{1}{c}{$\bolds{\eta}$} & \multicolumn{1}{c}{$\bolds{\alpha}$}
& \multicolumn{1}{c}{$\bolds{\tau^2}$} & \multicolumn{1}{c@{}}{$\bolds{\eta}$} \\
\hline
1 & 95.56 & $0.2526$ & $(-10.21, 10.40)$ & \hphantom{0}$(79.43, 119.54)$ &
$(0.2107, 0.2544)$ \\
2 & 156.90 & $0.1855$ & $(-3.19, 3.42)$ & $(125.96, 190.08)$ &
$(0.0922, 0.2257)$ \\
3 & 128.94 & $0.2476$ & $(-7.66, 7.76)$ & $(105.63, 159.54)$ &
$(0.1976, 0.2533)$ \\
4 & 129.92 & $0.2095$ & $(-3.57, 3.74)$ & $(104.54, 159.76)$ &
$(0.1264, 0.2380)$ \\
5 & 69.33 & $0.2522$ & $(-8.29, 8.23)$ & \hphantom{0}$(57.20, 86.44)$\hphantom{0} &
$(0.2091, 0.2543)$ \\
6 & 210.75 & $0.2542$ & $(-20.57, 19.68)$ & $(175.39, 268.45)$
& $(0.2136, 0.2549)$ \\
\hline
\end{tabular*}
\end{table}

GOF $p$-values resulting from the parametric bootstrap procedure of
Section \ref{sec6.3} are presented in Table \ref{table4} for each of
the four test statistics of Section~\ref{sec3.3}. Overall these values
provide no indication that we are able to dismiss model (\ref{eqnmsk})
as a plausible representation of the spatial structure present in these
data.

%
\begin{table}
\tablewidth=260pt
\caption{Parametric bootstrap $p$-values for the six agricultural field
trials} \label{table4}
\begin{tabular*}{\tablewidth}{@{\extracolsep{\fill}}lcccc@{}}
\hline
\textbf{Trial} & $\bolds{T_1}$ & $\bolds{T_2}$
& $\bolds{T_3}$ & $\bolds{T_4}$ \\
\hline
1 & 0.8348 & 0.7976 & 0.7086 & 0.7530 \\
2 & 0.3844 & 0.4182 & 0.2132 & 0.3262 \\
3 & 0.0852 & 0.1168 & 0.1506 & 0.1478 \\
4 & 0.1656 & 0.1084 & 0.1426 & 0.0972 \\
5 & 0.2162 & 0.1828 & 0.1754 & 0.2024 \\
6 & 0.3502 & 0.2382 & 0.4642 & 0.2984 \\
\hline
\end{tabular*}
\end{table}

\section{Conclusions}\label{sec7}

In this article we have introduced a practical method to assess the
aptness of Markov random field models for representing spatial
processes. This method is based on special sets of locations we have
called concliques that partition the total set of observed locations
such that generalized spatial residuals within each conclique
approximate realizations of independent random variables on the unit
interval. These generalized spatial residuals can be combined across
nonindependent concliques in natural ways to produce GOF statistics
that correspond to Gaussian empirical processes that have identifiable
limit distributions. While those limit distributions can involve
complex covariance structures, we have demonstrated that finite sample
versions of the GOF statistics appear to converge rather quickly to
their limits, at least in the case of a simple null hypothesis. This
implies that, in an application, approximation of their limit
distributions under a~suitable null hypothesis will provide a useful
reference distribution against which to compare the value of an
observed GOF statistic. The composite hypothesis setting introduces a
considerably more complicated situation than does the simple hypothesis
setting, because limit laws involve covariances that cannot be easily
determined either explicitly or numerically. In an application,
resampling methods would seem to hold the greatest promise for
approximating distributions of GOF statistics based on generalized
spatial residuals.
While developing spatial subsampling or block bootstraps (cf. \cite
{L03}, Chapter 12) for this purpose requires further investigation, the
use of such resampling was illustrated in this article in the
application to agricultural field trials.

We wish to comment on a number of issues that involve the distinction
between application of the GOF methodology developed and the production
of theoretical results for that methodology. First is the issue of
stationarity. There is nothing in the definition or construction of
generalized spatial residuals, or GOF statistics constructed from them,
that requires a stationary model. All that is needed is identification
of a full conditional distribution for each location (\ref{eqn1}) that
may then be used in (2.1), and assurance that a joint distribution
having these conditionals exists. Assumptions of stationarity made in
this article facilitate the production of theoretical results needed to
justify use of the methodology. Another issue is application to
discrete cases. While the data examples given have considered
continuous conditional models, we have applied random generalized
spatial residuals to models formed from Winsorized Poisson conditional
distributions \cite{KC97} with promising empirical results.
Similar to questions of stationarity and discrete cases, there is
nothing in the constructive methodology that requires a regular lattice
or that each location have the same number of neighbors. Use of
a~regular lattice in this article again facilitates the demonstration of
theoretical properties, but this is not needed to implement the
procedures suggested. The application of Section \ref{sec6} involved a regular
lattice, but no border strip or other boundary conditions were imposed
to render neighborhoods of equal size. It should certainly be
anticipated that there may be edge effects on GOF statistics as
developed here, just as there are edge effects on properties of
estimators. How severe these effects might be in various settings, and
whether the use of modified boundary conditions (e.g., \cite{BK95})
could mitigate such effects is an issue in need of additional
investigation. Essentially the same thoughts can be offered relative to
potential sparseness that might occur in an application. Locations
lacking neighbors entirely could be considered members of any conclique
one chooses, and construction of GOF statistics would proceed
unhindered. What the effects of varying degrees of sparseness are
remains an open question. Of course, if no locations have any
neighbors, then the methodology presented here reduces to statistics
constructed on the basis of the ordinary probability integral transform
for independent random variables.

As with all GOF tests, the procedure based on generalized spatial
residuals developed in this article serves as a vehicle for assessing a
selected model for overall adequacy, not as a vehicle for selection of
the most attractive model in the first place. This is important in
consideration of fitted models under the composite setting, in which we
can think of estimation as having ``optimized'' a given model structure
for description of a set of observed data. There may be two or more
such structures that could be, with the best choice of parameter values
possible, viewed as plausible data generating mechanisms for a set of
observations. This does not necessarily mean, however, that those
different structures are equally pleasing as models for the problem
under consideration.

Finally, we mention a connection with the assessment of $k$-step ahead
forecasts in a time series setting. Let $\{X_t; t \in\mathbb{Z}\}$
denote a series of random variables observed at discrete, equally
spaced, points in time. The probability integral transform with
distributions conditioned on the present and past has been used to
construct $k$-step ahead residuals $U_{t+k} = F_{t+k|t}(X_{t+k})$,
where the conditioning in $F$ is on $\{X_t, X_{t-1}, \ldots\}$ (e.g.,
\cite{Da84,Di98,Gn07}). While our use of the probability integral
transform is similar to what is done in this time series setting, the
conditioning requirements are quite distinct. In the spatial setting,
two spatial residuals are independent only if neither is in the
conditioning set of the other (i.e., are both in the same conclique).
In time series $k$-step ahead forecasts, two values, $U_i$ and $U_j$,
will be independent only if either $X_{i}$ \textit{is} in the conditioning
set of $X_{j}$, or vice versa. The difference stems from the use of
full conditionals in spatial Markov random field models, rather than the
sequential
conditionals in the time series context which need not invoke a Markov
property at all. The approach taken to the development of theoretical
results in this article could potentially be used in the time series
setting, but the modifications require further
investigation.\looseness=1

\section{Proof of generalized spatial residual properties}\label{sec8}

As Theorem \ref{Theorem2.1} provides the main distributional result for
generalized
spatial residuals (\ref{genresid}) from concliques, which are
fundamental to the GOF test statistics of Sections \ref{sec3}
and~\ref{sec4}, we
establish Theorem \ref{Theorem2.1} here. The proofs of other results
from the
manuscript are provided in supplementary material \cite{KLN11}.

Let $\mathcal{Q}$ denote a finite subset of
conclique $\C\subset\mathbb{Z}^d$ with
$|\mathcal{Q}|=l \geq2$, and
let $ \mathcal{I}_\mathcal{Q}= \bigcup_{\bfs\in\mathcal{Q}} \{
\bfi\dvtx
\bfi\in\mathcal{N}(\bfs)\} = \{\bfs_1,\ldots,\bfs_L\}$, $L\geq1$, be
the finite index set of all neighbors
of sites in $\mathcal{Q}$; additionally, enumerate the $l$ elements of
$\mathcal{Q}$ as $\mathcal{Q}=\{\bfs_{1+L},\ldots, \bfs_{l+L}\}$, say.
With respect to the enumeration of $ \mathcal{I}_\mathcal{Q}$ and
$\mathcal{Q}$, let $F_1(\cdot)$ denote the marginal c.d.f. of $Y(\mathbf
{s}_1)$, and let $F_j(\cdot)$, $2 \leq j \leq L+l$, denote the
conditional c.d.f. of $Y(\mathbf{s}_j)$ given $Y(\mathbf{s}_1),\ldots,
Y(\mathbf{s}_{j-1})$; define the function $F_j^{-}(\cdot)$
by the left limits of $F_j(\cdot)$. By the randomized PIT \cite{Br07},
$\{ (1-A(\bfs_j))\cdot F_j[Y(\mathbf{s}_j)] + A(\bfs_j)\cdot
F_j^{-}[Y(\mathbf{s}_j)]\dvtx j=1,\ldots,L+l\}$ are i.i.d. Uniform $(0,1)$
random variables.

For any $i,k \in\{L+1,\ldots,L+l\}$, variables $Y(\mathbf{s}_i)$ and
$Y(\mathbf{s}_k)$
belong to the conclique $\mathcal{Q}$ so that all neighbors of
$Y(\mathbf{s}_i)$
and $Y(\mathbf{s}_k)$ are among\vspace*{1pt} $\{Y(\mathbf{s}_j)\}_{j=1}^{L}$. By the
Markov property
({\ref{processmodel}), $F_j[Y(\mathbf{s}_j)] = F[Y(\mathbf{s}_j |\{
Y(\bfs)\dvtx\bfs\in\mathcal{N}(\mathbf{s}_j)\})]$ holds
and we may equivalently write (\ref{genresid}) as $U(\bfs_j)= (1-A(\bfs
_j))\cdot F_j[Y(\mathbf{s}_j)] + A(\bfs_j)\cdot F_j^{-}[Y(\mathbf{s}_j)]$
for any $j\in\{L+1,\ldots,L+l\}$, though these relationships may not
necessarily hold for $j=1,\ldots,L$.
Hence,
$\{U(\bfs)\dvtx\bfs\in\mathcal{Q}\}$ are i.i.d. Uniform $(0,1)$ variables
for any arbitrary finite subset $\mathcal{Q}$ of $\C$.

\section*{Acknowledgments}

The authors are very grateful to an Associate Editor and three referees
for thoughtful comments and suggestions which significantly clarified
and improved the manuscript.

\begin{supplement}[id=suppA]
\stitle{Proofs of main results for spatial GOF test statistics\\}
\slink[doi]{10.1214/11-AOS948SUPP} 
\sdatatype{.pdf}
\sfilename{aos948\_supp.pdf}
\sdescription{A supplement \cite{KLN11} provides proofs of all
asymptotic distributional results from Section \ref{sec4}, regarding the
conclique-based spatial GOF test statistics in simple and composite
null hypothesis settings (Proposition \ref{Proposition4.1}, Theorem
\ref{Theorem4.2}, Corollary \ref{Corollary4.3},
Theorem \ref{Theorem4.4}, Corollary \ref{Corollary4.5}). The proof in
the composite hypothesis
case is particularly nonstandard; see Section \ref{sec4.4}.}
\end{supplement}

%

\printaddresses

\end{document}